\documentclass{article}
\usepackage[utf8]{inputenc}
\usepackage[english]{babel}
\usepackage{amsthm}
\usepackage{amsmath}
\usepackage{graphicx}
\usepackage{url}
\usepackage{hyperref}

\begin{document}

\title{Prime product formulas for the Riemann zeta function and related identities}

\author{Artur Kawalec}

\date{}
\maketitle

\begin{abstract}
In this article, we derive a Euler prime product formula for the magnitude of the Riemann zeta function $\zeta(s)$ valid for $\Re(s)>1$, as well as similar formulas for $\zeta(s)$ valid for an even and odd $k$th positive integer argument. We shall further give a set of generated formulas for $\zeta(k)$ up to $11$th order, including  Ap\'ery's constant, and also construct formulas for $\zeta(3/2)$. We'll also validate these formulas numerically.

\end{abstract}

\section{Main Prime Product Formula}
The Euler's prime product formula is a key connection between the Riemann zeta function and prime numbers. If $p_n$ is a sequence of $n$th prime numbers denoted such that $p_1=2$, $p_2 = 3$, $p_3=5$ and so on, then the Riemann zeta function is given as Euler prime product
\begin{equation}\label{eq:1}
\zeta(s) = \prod_{n=1}^\infty\ \left(1-\frac{1}{p^s_n}\right)^{-1},
\end{equation}
which converges absolutely for $\Re(s)>1$.  Next, we shall seek to evaluate the magnitude $\mid \zeta(s) \mid$ for complex argument $s=\sigma+it$ for which $\sigma >1$. First, by substituting complex argument we have
\begin{equation}\label{eq:2}
\zeta(\sigma+it) = \prod_{n=1}^\infty\ \left(1-\frac{1}{p^{\sigma+it}_n}\right)^{-1},
\end{equation}
and by further algebraic simplification we arrive at
\begin{equation}\label{eq:3}
\zeta(\sigma+it) = \prod_{n=1}^\infty\ \frac{p_n^{\sigma}-e^{it\log{p_n}}}{p_n^{\sigma}+p_n^{-\sigma}-2\cos{(t\log{p_n}})}.
\end{equation}
The magnitude can then be written as
\begin{equation}\label{eq:4}
\mid \zeta(\sigma+it) \mid^2 = \prod_{n=1}^\infty (1+p_n^{-2\sigma})^{-1}\left(1-\frac{2}{p_n^{\sigma}+p_n^{-\sigma}}\cos(t\log p_n)\right)^{-1}.
\end{equation}
Furthermore, using the identity

\begin{equation}\label{eq:5}
\frac{\zeta(2\sigma)}{\zeta(\sigma)}= \prod_{n=1}^\infty (1+p_n^{-\sigma})^{-1}, 
\end{equation}
we obtain the main formula:

\begin{equation}\label{eq:6}
\mid \zeta(\sigma+it) \mid^2 = \frac{\zeta(4\sigma)}{\zeta(2\sigma)}\prod_{n=1}^\infty \left(1-\frac{2}{p_n^{\sigma}+p_n^{-\sigma}}\cos(t\log p_n)\right)^{-1}.
\end{equation}
Alternatively, the prime factor terms can also be expressed as hyperbolic cosines, thus resulting in

\begin{equation}\label{eq:7}
\mid \zeta(\sigma+it) \mid^2 = \frac{\zeta(4\sigma)}{\zeta(2\sigma)}\prod_{n=1}^\infty \left(1-\frac{\cos(t\log p_n)}{\cosh(\sigma\log p_n)}\right)^{-1}.
\end{equation}
This completes the derivation of the magnitude of $\zeta(\sigma+it)$. In the next section, we will use these results and derive an integer formula for the Riemann zeta function.

\section{Prime Product Integer Formula}
Using the result of the previous section, we let $t=0$ in equation (6) formulation which results in

\begin{equation}\label{eq:8}
\zeta(\sigma)^2 = \frac{\zeta(4\sigma)}{\zeta(2\sigma)}\prod_{n=1}^\infty \left(1-\frac{2}{p_n^{\sigma}+p_n^{-\sigma}}\right)^{-1},
\end{equation}
which is valid for any $\sigma>1$. If for a positive integer $k$ we let $\sigma=k$ and using the well-known identity for

\begin{equation}\label{eq:9}
\zeta(2k) = \frac{(-1)^{k+1}B_{2k}(2\pi)^{2k}}{2(2k)!},
\end{equation}
where $B_k$ are Bernoulli numbers, then the zeta terms simplify as

\begin{equation}\label{eq:10}
\frac{\zeta(4k)}{\zeta(2k)}=\pi^{2k}{\frac{(-1)^{k}2^{2k}B_{4k}(2k)!}{B_{2k}(4k)!}}, 
\end{equation}
and hence, the integer formula is obtained as

\begin{equation}\label{eq:11}
\zeta(k) =\pi^{k}\sqrt{\frac{(-1)^{k}2^{2k}B_{4k}(2k)!}{B_{2k}(4k)!}}\prod_{n=1}^\infty \left(1-\frac{2}{p_n^{k}+p_n^{-k}}\right)^{-1/2}.
\end{equation}
Alternatively, we also have another form:

\begin{equation}\label{eq:12}
\zeta(k) =\pi^{k}\sqrt{\frac{(-1)^{k}2^{2k}B_{4k}(2k)!}{B_{2k}(4k)!}}\prod_{n=1}^\infty\frac{\sqrt{p_n^{2k}+1}}{p_n^k-1}
\end{equation}
by simplifying the prime terms. But if one wishes to evaluate the complex magnitude, then the $\cos(t\log(p_n))$ term must be included. These formulas can be used to evaluate $\zeta(k)$ for any positive integer greater than unity as will be presented in the next section.

\section{Evaluation of $\zeta(2)$}
By setting $k=2$ into the above formula yields

\begin{equation}\label{eq:13}
\zeta(2) = \frac{\pi^2}{\sqrt{105}}\prod_{n=1}^\infty \left(1-\frac{2}{p_n^{2}+p_n^{-2}}\right)^{-1/2}.
\end{equation}
The numerical computation of this formula converges to the correct result. We also note that the prime product term converges to a constant factor

\begin{equation}\label{eq:14}
\prod_{n=1}^\infty \left(1-\frac{2}{p_n^{2}+p_n^{-2}}\right)^{-1/2} \to \frac{\sqrt{105}}{6},
\end{equation}
where we obtain the familiar result
\begin{equation}\label{eq:15}
\zeta(2) = \frac{\pi^2}{6}.
\nonumber
\end{equation}
In the next example, we evaluate the magnitude of $\zeta(2+i)$ using equation (7) as
\begin{equation}\label{eq:16}
\mid \zeta(2+i) \mid = \frac{\pi^2}{\sqrt{105}}\prod_{n=1}^\infty \left(1-\frac{\cos(\log p_n)}{\cosh(2\log p_n)}\right)^{-1/2} \to 1.23075241321861\dots .
\nonumber
\end{equation}

\section{Evaluation of $\zeta(3)$}
Perhaps of greater interest is the formula for $\zeta(3)$, or Ap\'ery's constant, for which there is no known representation as a rational multiple of $\pi^3$ such as for an even order case $\zeta(2k)$ by equation (9). By setting $k=3$ into the above formula directly results in

\begin{equation}\label{eq:17}
\zeta(3) = \pi^3 \sqrt{\frac{691}{675675}}\prod_{n=1}^\infty \left(1-\frac{2}{p_n^{3}+p_n^{-3}}\right)^{-1/2}.
\end{equation}
Numerical computation also validate the convergence to correct a value. We note that the prime product term converges to a constant factor
\begin{equation}\label{eq:18}
\prod_{n=1}^\infty \left(1-\frac{2}{p_n^{3}+p_n^{-3}}\right)^{-1/2} \to 1.21228661439701\ldots,
\nonumber
\end{equation}
which is approximately equal to $\zeta(3)$ to within one decimal place, but it is not known if it can be expressed in terms of known constants. As another example, we compute the magnitude of $\zeta(3+i)$ as
\begin{equation}\label{eq:19}
\mid \zeta(3+i) \mid = \pi^3 \sqrt{\frac{691}{675675}}\prod_{n=1}^\infty \left(1-\frac{\cos(\log p_n)}{\cosh(3\log p_n)}\right)^{-1/2} \to 1.11710067922572\dots.
\nonumber
\end{equation}
In Appendix A, we summarized these formulas up to $11$th order. And in Appendix B, we summarized the numerical validation of these formulas in Table 1 up to $15$ decimal places.

\section{Evaluation of $\zeta(3/2)$}

We can use the same approach to obtain a formula for $\zeta(3/2)$, by substituting $\sigma=3/2$ into equation (11) yields

\begin{equation}\label{eq:1}
\zeta(3/2)^2 = \frac{\zeta(6)}{\zeta(3)}\prod_{n=1}^\infty \left(1-\frac{2}{p_n^{3/2}+p_n^{-3/2}}\right)^{-1},
\end{equation}
and using $\zeta(3)$ formula obtained earlier results in
\begin{equation}\label{eq:20}
\zeta(3/2) = \pi^{3/2} \sqrt[4]{\frac{675675}{617080275}}\prod_{n=1}^\infty \left(1-\frac{2}{p_n^{3/2}+p_n^{-3/2}}\right)^{-1/2} \left(1-\frac{2}{p_n^{3}+p_n^{-3}}\right)^{1/4}
\end{equation}
with leading $\pi^{3/2}$ factor. Numerical validation of this formula is also summarized in Appendix B. And similarly, as before, we evaluate the magnitude of $\zeta(3/2+i)$ as
\begin{equation}\label{eq:21}
\mid \zeta(3/2+i) \mid = \frac{\pi^3}{3\sqrt{105\zeta(3)}}\prod_{n=1}^\infty \left(1-\frac{\cos(\log p_n)}{\cosh(3/2\log p_n)}\right)^{-1/2} \to  1.2536382542739\dots,
\nonumber
\end{equation}
but without expanding $\zeta(3)$, which is now in terms of $\pi^{3}$ constant.

\section{Other Prime Product Formulas}

We obtain another set of similar formulas as such. Using the identity (5), and multiplying both sides by the Euler prime product equation (1) results in

\begin{equation}\label{eq:22}
\zeta(s)^2 = \zeta(2s)\prod_{n=1}^\infty (1+p_n^{-s})(1-p_n^{-s})^{-1}
\end{equation}
and so this leads to
\begin{equation}\label{eq:23}
\zeta(s)^2 = \zeta(2s)\prod_{n=1}^\infty \left(\frac{p_n^s+1}{p_n^s-1}\right),
\end{equation}
and if $s$ is a positive integer $k$, then
\begin{equation}\label{eq:24}
\zeta(k) = \pi^k \sqrt{\frac{(-1)^{k+1}B_{2k}2^{2k-1}}{(2k)!}} \prod_{n=1}^\infty \left(\frac{p_n^k+1}{p_n^k-1}\right)^{1/2}.
\end{equation}
This results in
\begin{equation}\label{eq:25}
\zeta(2) = \frac{\pi^2}{\sqrt{90}} \prod_{n=1}^\infty \left(\frac{p_n^2+1}{p_n^2-1}\right)^{1/2}
\end{equation}
and similarly, for Ap\'ery's constant we have
\begin{equation}\label{eq:26}
\zeta(3) = \frac{\pi^3}{3\sqrt{105}} \prod_{n=1}^\infty \left(\frac{p_n^3+1}{p_n^3-1}\right)^{1/2}.
\end{equation}
And just as before, we obtain an identity for $\zeta(3/2)$ as

\begin{equation}\label{eq:27}
\zeta(3/2) = \sqrt{\zeta(3)}\prod_{n=1}^\infty \left(\frac{p_n^{3/2}+1}{p_n^{3/2}-1}\right)^{1/2}
\end{equation}
and
\begin{equation}\label{eq:28}
\zeta(3/2) = \frac{\pi^{3/2}}{\sqrt[4]{945}} \prod_{n=1}^\infty \left(\frac{p_n^{3/2}+1}{p_n^{3/2}-1}\right)^{1/2} \left(\frac{p_n^3+1}{p_n^3-1}\right)^{1/4}.
\end{equation}
Higher-order formulas easily follow.  Although this form cannot be used to compute the magnitude, such as by equations (6) and (7), it will, however, simplify to the original Euler product (1) if one substitutes back the $\zeta(2k)$ by equation (9).

\section{Conclusion}
A simple formula for the magnitude of the Riemann zeta function was derived based on Euler prime products, which result in a variety of formulas, such as for positive integer argument $k>1$. We also notice that for an even order case, the prime product term in equation (11) simplifies to the standard $\zeta(2k)$ result given by equation (9), such as $\pi^2/6$ for $k=2$. For an odd order case, however, it is not known if the prime product term can be expressed in closed-form, such as in the $\zeta(2k)$ case. The usefulness of these formulas is that the magnitude of $\zeta(s)$ for complex argument with $\Re(s)>1$ can be evaluated by using the $\cos(t\log p_n)$ term in equations (6)(7) and (11). We also derived a similar set of formulas from equation (5), which are obtained by multiplying the Euler prime product by $\prod_{n=1}^\infty (1+p_n^s)^{-1}$. The main theme behind these formulas is that $\pi^k$ term is extracted from the Euler prime product formula and a closed-form representation of $\zeta(2k)$, and by combining multiple formulas many different varieties can be obtained.
\newpage
\section{Appendix A}
The evaluation of equation (11) in the Mathematica software package up to the $11$th order:

\begin{equation}\label{eq:29}
\zeta(2) = \frac{\pi^2}{\sqrt{105}}\prod_{n=1}^\infty \left(1-\frac{2}{p_n^{2}+p_n^{-2}}\right)^{-1/2}
\nonumber
\end{equation}

\begin{equation}\label{eq:30}
\zeta(3) = \pi^3 \sqrt{\frac{691}{675675}}\prod_{n=1}^\infty \left(1-\frac{2}{p_n^{3}+p_n^{-3}}\right)^{-1/2}
\nonumber
\end{equation}

\begin{equation}\label{eq:31}
\zeta(4) = \pi^4 \sqrt{\frac{3617}{34459425}}\prod_{n=1}^\infty \left(1-\frac{2}{p_n^{4}+p_n^{-4}}\right)^{-1/2}
\nonumber
\end{equation}

\begin{equation}\label{eq:32}
\zeta(5) = \pi^5 \sqrt{\frac{174611}{16368226875}}\prod_{n=1}^\infty \left(1-\frac{2}{p_n^{5}+p_n^{-5}}\right)^{-1/2}
\nonumber
\end{equation}

\begin{equation}\label{eq:33}
\zeta(6) = \pi^6 \sqrt{\frac{236364091}{218517792968475}}\prod_{n=1}^\infty \left(1-\frac{2}{p_n^{6}+p_n^{-6}}\right)^{-1/2}
\nonumber
\end{equation}

\begin{equation}\label{eq:34}
\zeta(7) = \pi^7 \sqrt{\frac{3392780147}{30951416768146875}}\prod_{n=1}^\infty \left(1-\frac{2}{p_n^{7}+p_n^{-7}}\right)^{-1/2}
\nonumber
\end{equation}

\begin{equation}\label{eq:35}
\zeta(8) = \pi^8 \sqrt{\frac{7709321041217}{694097901592400930625}}\prod_{n=1}^\infty \left(1-\frac{2}{p_n^{8}+p_n^{-8}}\right)^{-1/2}
\nonumber
\end{equation}

\begin{equation}\label{eq:36}
\zeta(9) = \pi^9 \sqrt{\frac{26315271553053477373}{23383376494609715287281703125}}\prod_{n=1}^\infty \left(1-\frac{2}{p_n^{9}+p_n^{-9}}\right)^{-1/2}
\nonumber
\end{equation}

\begin{equation}\label{eq:37}
\zeta(10) = \pi^{10} \sqrt{\frac{261082718496449122051}{2289686345687357378035370971875}}\prod_{n=1}^\infty \left(1-\frac{2}{p_n^{10}+p_n^{-10}}\right)^{-1/2}
\nonumber
\end{equation}

\begin{equation}\label{eq:38}
\zeta(11) = \pi^{11} \sqrt{\frac{2530297234481911294093}{219012470258383844016431785453125}}\prod_{n=1}^\infty \left(1-\frac{2}{p_n^{11}+p_n^{-11}}\right)^{-1/2}.
\nonumber
\end{equation}

\newpage
\section{Appendix B}

In the following table, we summarize the numerical computation of equation (11) in the Mathematica software package to $15$ decimal places for the first $1000$ prime product terms. We also note that convergence is slower for smaller arguments, such as $k=1.5$ or $k=2$,  and that all the higher-order argument converged faster.

\begin{table}[ht]
\caption{Evaluation of $\zeta(k)$ and the new formula for $\zeta(k)$ respectively} 
\centering 
\begin{tabular}{c c c} 
\hline\hline 
k & $\zeta(k)$ & $\zeta(k)$ Equation (11) \\ [0.5ex] 
\hline 
1.5  & 2.61237534868549 & 2.60691093229650 \\ 
2  & 1.64493406684823 & 1.64491317470628 \\
3  & 1.20205690315959  & 1.20205690215259 \\
4  & 1.08232323371114   & 1.08232323371106 \\
5  & 1.03692775514337  & 1.03692775514337 \\
6  & 1.01734306198445 & 1.01734306198445 \\
7  & 1.00834927738192 & 1.00834927738192 \\
8  & 1.00407735619794  & 1.00407735619794 \\
9  &  1.00200839282608  & 1.00200839282608 \\
10  & 1.00099457512782  & 1.00099457512782 \\
11 & 1.00049418860412  & 1.00049418860412
 \\ [1ex] 
\hline 
\end{tabular}
\label{table:nonlin} 
\end{table}

\texttt{Email: art.kawalec@gmail.com}

\begin{thebibliography}{9}
\bibitem{latexcompanion}
H.M. Edwards.
\textit{Riemann's Zeta Function}.
Dover Publication, Mineola, New York 1974.

\end{thebibliography}
\end{document}